\documentclass{article}
\pdfoutput=1
\usepackage{amssymb}	
\usepackage{amsmath}    
\usepackage{setspace} 
\usepackage{latexsym}	
\usepackage{longtable}
\usepackage{multirow}
\usepackage{epsfig}
\allowdisplaybreaks

\newtheorem{algorithm}{Algorithm}[section]

\newtheorem{theorem}{Theorem}[section]

\newtheorem{remark}{Remark}[section]

\newcommand{\qed}{\nobreak \ifvmode \relax \else \ifdim\lastskip<1.5em \hskip-\lastskip \hskip1.5em plus0em minus0.5em \fi \nobreak \vrule height0.75em width0.5em depth0.25em\fi} 

\def\A{{\bf A}}
\def\B{{\bf B}}
\def\C{{\bf C}}

\def\F{{\bf F}}
\def\G{{\bf G}}
\def\H{{\bf H}}
\def\I{{\bf I}}
\def\J{{\bf J}}
\def\K{{\bf K}}

\def\0{{\bf 0}}

\def\a{{\bf a}}
\def\b{{\bf b}}

\def\e{{\bf e}}

\def\g{{\bf g}}
\def\h{{\bf h}}

\def\q{{\bf q}}

\def\s{{\bf s}}
\def\t{{\bf t}}
\def\u{{\bf u}}

\def\x{{\bf x}}

\def\Tr{{\rm T}}

\def\Re{{\cal R}_e}

\def\diag{{\rm diag}}

\title{Singularity-Free Spacecraft Attitude Control Using Variable-Speed Control Moment Gyroscopes}
\author {Yaguang Yang\thanks{
Office of Research, NRC, 21 Church Street, Rockville, 20850. Email:
yaguang.yang@verizon.net} 
}
\date{\today}

\begin{document}

\maketitle  

\begin{abstract}
This paper discusses spacecraft control using variable-speed CMGs.
A new operational concept for VSCMGs is proposed. This new concept
makes it possible to approximate the complex nonlinear system by
a linear time-varying system (LTV). As a result, an effective control 
system design method, Model Predictive Control (MPC) using robust
pole assignment, can be used to design the spacecraft
control system using VSCMGs. A nice feature of
this design is that the control system does not have any singular 
point. A design example is provided. The simulation result shows
the effectiveness of the proposed method.
\end{abstract}


\newpage

\section{ Introduction}

Control Moment Gyros (CMGs) are an important type of actuators
used in spacecraft control because of their well-known torque
amplification property \cite{kurokawa98}. The conventional use of
CMG keeps the flywheel spinning in a constant speed, while torques of the CMG are produced by changing the gimbal's 
rotational speed \cite{jt04}. A more complicated operational 
concept is the so-called variable-speed control gimbal gyros 
(VSCMG) in which the flywheel's speed of the CMG is allowed to be 
changed too. This idea was first proposed by Ford in his Ph.D 
dissertation \cite{ford97} where he derived a mathematical model 
for VSCMGs which is now widely used in literatures. Because of 
the extra freedom of VSCMG, it can generate torques on a plane
perpendicular to the gimbal axis while the conventional CMG can
only generate a torque in a single direction at any instant of
time \cite{yt06}.

The existing designs of spacecraft control system using CMG or 
VSCMG rely on the calculation of the desired torques and 
then determines the VSCMG's gimbal speed and flywheel speed. 
This designs have a fundamental problem because there are singular 
points where the gimbal speed and flywheel speed cannot be found
given the desired torques. Extensive literatures focus on this
difficulty of implementation in the last few decades, for example,
\cite{kurokawa98,fh00,yt04,kurokawa07,zmmt15} and references
therein. Another difficulty associated with the control system
design using CMG or VSCMG is that the nonlinear dynamical models 
for these type of actuators are much more complicated than other
types of actuators used for spacecraft attitude control systems. 
Most proposed designs, for example \cite{jt04,yt06,fh00,yt02,
fh97,svj98,ma13,jh11}, use Lyapunov stability theory 
for nonlinear systems. There are two shortcomings of this
design method: first, there is no systematic way to find the
desired Lyapunov function, and second, the design does not 
consider the system performance but only stability.

In this paper, we propose a different operational concept for
VSCMG: the flywheels of the cluster of the VSCMG do not always
spin at high speed, they spin at high speed only when they need
to. The same is true for the gimbals. This operational strategy 
makes the origin (the state variables at zero) an equilibrium 
point, where a linearized model can be established. Therefore, 
some mature linear system design methods can be used and system 
performance can be part of the design by using these linear system 
design methods. Additional advantages of the proposed operational 
concept are: (a) energy saving due to normally reduced spin speed 
of flywheels and gimbals therefore reduction of operational cost, 
(b) seamless implementation (singularity free) because the 
control of the spacecraft is achieved by accelerating or 
decelerating the flywheels and gimbals, therefore, there is no 
inverse from desired torques to the speeds of the gimbals and
flywheels.

It is worthwhile to point out that the linearized model is a 
linear time-varying (LTV) system. The design methods for linear 
time-invariant (LTI) systems cannot be directly applied to LTV
systems. A popular design method for LTV system is 
the so-called gain scheduling design method, which has been
discussed in several decades, for example, \cite{rugh90,
rs00,lr90,sb92}. The basic idea is to fix the time-varying
model in a number of ``frozen'' models and using linear system
design method for each of these ``frozen'' linear time-invariant 
systems. When the parameters of the LTV system are not in these
``frozen'' points, interpolation is used to calculate the
feedback gain matrix. 

Although, gain scheduling design has been proved to be
effective for many applications for LTV systems, it has an
intrinsic limitation for some time-varying systems which have
many independent time-varying variables, which is the
case for spacecraft control using VSCMGs. As we will see
later that this control system matrices $(\A,\B)$ have
many independent time-varying parameters and the 
computation of the gain scheduling design is too much to be
feasible. Therefore, we will consider another popular
control system design method, the so-called Model
Predictive Control (MPC) \cite{aw13}. According to
a theorem in \cite{rugh93}, under certain conditions, the 
closed-loop LTV system designed by MPC method is stable. To 
meet some of the required stability conditions imposed on 
the LTV system \cite{rugh93}, we propose using the robust 
pole assignment design \cite{yt93,tity96} for the MPC design.

The remainder of the paper is organized as follows. Section 2 
derives the spacecraft model using variable-speed CMG. Section
3 discusses both gain scheduling design and the MPC design method 
for spacecraft control using variable-speed CMG. This analysis provides a technical basis of selecting the MPC design over gain
scheduling design for this problem. Section 4 provides a design 
example and simulation result. Section 5 is the summary of the conclusions of the paper.

\section{Spacecraft model using variable-speed CMG} 

Throughout the paper, we will repeatedly use a skew-symmetric 
matrix which is related to the cross product of two vectors. Let 
$\a=[ a_1,a_2,a_3 ]^{\Tr}$ and $\b=[ b_1,b_2,b_3 ]^{\Tr}$ be
two three dimensional vectors. We denote a matrix
\[
\a^{\times} =\left[ \begin{array}{ccc}
0 & -a_3 & a_2 \\ a_3 & 0 & -a_1 \\ -a_2 & a_1 & 0
\end{array} \right]
\]
such that the cross product of $\a$ and $\b$ is equivalent to a matrix
and vector multiplication, i.e., $\a \times \b = \a^{\times} \b$. 
  
\begin{figure}[ht]
\centering
\includegraphics[height=8cm,width=8cm]{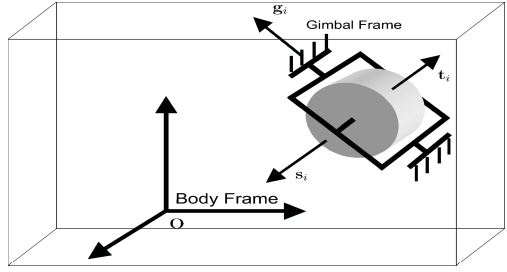}
\caption{Spacecraft body with a single VSCMG.}
\label{fig:gimbalFrame}
\end{figure}

Assuming that there are $N$ variable-speed CMGs installed in a 
spacecraft, following the notations of \cite{ford97}, we define 
a matrix 
\begin{equation}
\A_s = [ \s_{1},  \s_{2}, \ldots,  \s_{N}]
\end{equation}
such that the columns of $\A_s$, $\s_{j}$ ($j=1, \ldots, N$), 
specify the unit spin axes of the wheels in the spacecraft body 
frame. Similarly, we define 
$\A_g= [ \g_{1},  \g_{2}, \ldots,  \g_{N}]$ the matrix whose 
columns are the unit gimbal axes and 
$\A_t= [ \t_{1},  \t_{2}, \ldots,  \t_{N}]$ 
the matrix whose columns are the unit axes of the transverse 
(torque) directions, both are represented in the spacecraft 
body frame. Whereas $\A_g$ is a constant matrix, the matrices 
$\A_s$ and $\A_t$ depend on the gimbal angles. Let 
$\boldsymbol{\gamma}=[\gamma_1, \ldots, \gamma_N]^{\Tr}
\in [0,2\pi] \times \cdots  \times [0,2\pi]:=\Pi$ 
be the vector of $N$ gimbal angles,
\begin{equation}
[ \dot{\gamma}_1, \ldots,  \dot{\gamma_N}]^{\Tr}
= \dot{\boldsymbol{\gamma}}:=\boldsymbol{\omega}_g
= [ {\omega}_{g_1}, \ldots,  {\omega}_{g_N}]^{\Tr}
\label{dotGamma}
\end{equation}
be the vector of $N$ gimbal speed, then the following relations 
hold \cite{yt04} (see Figure \ref{fig:gimbalFrame}).


\begin{equation}
\dot{\s}_i = \dot{{\gamma}}_i \t_i={{\omega}}_{g_i} \t_i,  \hspace{0.1in}
\dot{\t}_i = -\dot{{\gamma}}_i \s_i={{\omega}}_{g_i} \s_i,  \hspace{0.1in}
\dot{\g}_i = 0.
\label{CMGcord}
\end{equation}
Denote
\begin{equation}
\boldsymbol{\Gamma}^c = \diag (\cos(\boldsymbol{\gamma})),  \hspace{0.1in}
\boldsymbol{\Gamma}^s = \diag (\sin(\boldsymbol{\gamma})).
\label{diagcs}
\end{equation}
A different but related expression is given in \cite{ford97}
\footnote{There are some typos in the signs in \cite{ford97} which 
are corrected in (\ref{ford1}) and (\ref{CMGmatrixDerivative}).}. 
Let $\A_{s_0}$ and $\A_{t_0}$ be initial spin axes and gimbal 
axes matrices at $\boldsymbol{\gamma_0}=\0$, then
\begin{subequations}
\begin{align}
\A_s(\boldsymbol{\gamma}) = \A_{s_0} \boldsymbol{\Gamma}^c+\A_{t_0} \boldsymbol{\Gamma}^s,  \\
\A_t(\boldsymbol{\gamma}) = \A_{t_0} \boldsymbol{\Gamma}^c-\A_{s_0} \boldsymbol{\Gamma}^s.
\end{align}
\label{ford1}
\end{subequations}
This gives
\begin{subequations}
\begin{align}
\dot{\A}_s = \A_{t_0} \diag (\dot{\boldsymbol{\gamma}})
= \A_{t_0} \diag ( \boldsymbol{\omega}_g),  \\
\dot{\A}_t = -\A_{s_0} \diag (\dot{\boldsymbol{\gamma}}) 
= -\A_{s_0} \diag (\boldsymbol{\omega}_g),
\end{align}
\label{CMGmatrixDerivative}
\end{subequations}
which are identical to the formulas of (\ref{CMGcord}).
Let $J_{s_j}$, $J_{g_j}$, and $J_{t_j}$ be the spin axis inertia, the gimbal
axis inertia, and the transverse axis inertia of the $j$-th CMG, let three 
$N \times N$ matrices be defined as 
\begin{equation}
\J_s = \diag (J_{s_j}),  \hspace{0.1in}
\J_g = \diag (J_{g_j}),  \hspace{0.1in}
\J_t = \diag (J_{t_j}).
\end{equation}
Let the spacecraft inertia matrix be $\J_b$, then the total inertia matrix
including the CMG clusters is given by \cite{ford97}
\begin{equation}
\J =\J_b +\A_s \J_s \A_s^{\Tr}+\A_g \J_g \A_g^{\Tr}+\A_t \J_t \A_t^{\Tr},
\end{equation}
which is a function of $\boldsymbol{\gamma}$, 
and $\boldsymbol{\gamma}$ is a variable depending of time $t$. 
Therefore, $\J$ is an implicit function of 
time $t$. Although $\J$ is a function $\boldsymbol{\gamma}$, 
the dependence of $\J$ on $\boldsymbol{\gamma}$ 
is weak, especially when the size of spacecraft main body is large \cite{yt06}.
We therefore assume that $\dot{\J}=0$ as treated in
\cite{yt02,fh97,svj98,ma13}. 
Let $\boldsymbol{\omega}=[ \omega_1,\omega_2,  \omega_3]^{\Tr}$ 
be the spacecraft body angular rate with respect to the inertial 
frame, $\boldsymbol{\beta}=[\beta_1, \ldots, \beta_N]^{\Tr}$ 
be the vector of $N$ wheel angles,
\begin{equation}
[ \dot{\beta}_1, \ldots,  \dot{\beta}_{N}]^{\Tr}
= \dot{ \boldsymbol{\beta}} :=\boldsymbol{\omega}_s
= [ {\omega}_{s_1}, \ldots,  {\omega}_{s_N}]^{\Tr}
\end{equation}
be the vector of $N$ wheel speed.
Denote 
\begin{equation}
\h_s=[J_{s_1} \dot{\beta}_{1} , \ldots, J_{s_N} \dot{\beta}_{N} ]^{\Tr}
=\J_s \boldsymbol{\omega}_s,
\end{equation}
\begin{equation}
\h_g=[J_{g_1} \dot{\gamma}_{1},\ldots, J_{g_N} \dot{\gamma}_{N}]^{\Tr}
=\J_g \boldsymbol{\omega}_g,
\end{equation}
and $\h_t$ be the $N$ dimensional vectors representing the 
components of absolute angular momentum of the CMGs about their 
spin axes, gimbal axes, and transverse axes respectively. 
The total angular momentum of the spacecraft with a cluster 
of CMGs represented in the body frame is given as
\begin{equation}
\h = \J_b \boldsymbol{\omega} +\sum_{i=1}^{N} 
\s_i J_{s_i} \dot{\beta}_{i} +\sum_{i=1}^{N} 
\g_i J_{g_i} \dot{\gamma}_{i}
= \J_b \boldsymbol{\omega} +\A_s \h_s + \A_g \h_g
= \J_b \boldsymbol{\omega} +\A_s \J_s \boldsymbol{\omega}_s
+ \A_g \J_g \boldsymbol{\omega}_g.
\label{CMGh}
\end{equation}
Taking derivative of (\ref{CMGh}) and using (\ref{CMGcord}) and 
$\dot{\J}=0$, noticing that gimbal axes are fixed, we have
\begin{eqnarray}
\dot{\h} & = & \J_b \dot{\boldsymbol{\omega} }
 +\sum_{i=1}^{N} \left( \dot{\s}_i J_{s_i} \dot{\beta}_{i} 
+ {\s_i} J_{s_i} \ddot{\beta}_{i} \right)
+\sum_{i=1}^{N} \left( \dot{\g}_i J_{g_i} \dot{\gamma}_{i}
+ {\g}_i J_{g_i} \ddot{\gamma}_{i} \right) \nonumber \\
 & = & \J_b \dot{\boldsymbol{\omega} }
 +\sum_{i=1}^{N}  \left(  \dot{\gamma}_i \t_i  J_{s_i} \dot{\beta}_{i} 
+ {\s_i} J_{s_i} \ddot{\beta}_{i} \right)
+\sum_{i=1}^{N} {\g}_i J_{g_i} \ddot{\gamma}_{i}  \nonumber \\
 & = & -{\boldsymbol{\omega}} \times \h + \t_e,
\end{eqnarray}
where $\t_e$ is the external torque. Denote 
$\boldsymbol{\Omega}_s=\diag(\boldsymbol{\omega}_s)$ and
$\boldsymbol{\Omega}_g=\diag(\boldsymbol{\omega}_g)$.
This equation can be written as a compact form as follows.
\begin{eqnarray}
\J_b \dot{\boldsymbol{\omega} }
 +   \A_t  \J_{s} \boldsymbol{\Omega}_{s} \boldsymbol{\omega}_g
+ \A_s \J_{s} \dot{\boldsymbol{\omega}}_{s} 
+ {\A}_g \J_{g} \dot{\boldsymbol{\omega}}_{g}  
= -{\boldsymbol{\omega}} \times 
(\J_b \boldsymbol{\omega} +\A_s \J_s \boldsymbol{\omega}_s
+ \A_g \J_g \boldsymbol{\omega}_g) + \t_e,
\end{eqnarray}
Note that the torques generated by wheel acceleration or
deceleration in the directions defined by $\A_s$ are given by
\begin{equation}
\t_s=-\J_s \boldsymbol{\dot{\omega}}_s 
= [ t_{s_1}, \ldots, t_{s_N} ]^{\Tr}
\label{CMGtw}
\end{equation}
(note that vectors $\t_i$ in $\A_t$ are axes and scalars $t_{s_i}$
in $\t_s$ are torques) and the torques generated by gimbal acceleration or deceleration in the directions defined by $\A_g$ 
are given by
\begin{equation}
\t_g=-\J_g \boldsymbol{\dot{\omega}}_g
= [ t_{g_1}, \ldots, t_{g_N} ]^{\Tr},
\label{CMGtg}
\end{equation}
the dynamical equation can be expressed as
\begin{equation}
\J_b \dot{\boldsymbol{\omega} }
 +   \A_t  \J_{s} \boldsymbol{\Omega}_{s} \boldsymbol{\omega}_g
 + {\boldsymbol{\omega}} \times 
(\J_b \boldsymbol{\omega} +\A_s \J_s \boldsymbol{\omega}_s
+ \A_g \J_g \boldsymbol{\omega}_g)
=  \A_s \t_s  + {\A}_g \t_{g}+ \t_e.
\label{CMGdynamics}
\end{equation}
Let 
\begin{equation}
\bar{\q}=[q_0, q_1, q_2, q_3]^{\Tr}=[q_0, \q^{\Tr}]^{\Tr}=
\left[ \cos(\frac{\alpha}{2}), \hat{\e}^{\Tr}\sin(\frac{\alpha}{2}) \right]^{\Tr}
\end{equation}
be the quaternion representing the rotation of the body frame relative to the inertial frame, 
where $\hat{\e}$ is the unit length rotational axis and $\alpha$ is the rotation angle about $\hat{\e}$. 
Therefore, the reduced kinematics equation becomes \cite{yang10}
\begin{eqnarray} \nonumber
\left[  \begin{array} {c} \dot{q}_1 \\ \dot{q}_2 \\ \dot{q}_3
\end{array} \right]
& = & \frac{1}{2}  \left[  \begin{array} {ccc} 
\sqrt{1-q_1^2-q_2^2-q_3^2} & -q_3  & q_2 \\
q_3 & \sqrt{1-q_1^2-q_2^2-q_3^2} & -q_1 \\
-q_2  &  q_1  & \sqrt{1-q_1^2-q_2^2-q_3^2} \\
\end{array} \right] 
\left[  \begin{array} {c}  \omega_{1} \\ \omega_{2} \\ \omega_{3}
\end{array} \right]  \\
& = & \g(q_1,q_2, q_3, \boldsymbol{\omega}),
\label{nadirModel2}
\end{eqnarray}
or simply
\begin{equation}
\dot{\q}=\g(\q, \boldsymbol{\omega}).
\label{gfunction}
\end{equation}
The nonlinear time-varying spacecraft control system model 
can be written as follows:
\begin{eqnarray}
\left[ \begin{array}{c}
\dot{\boldsymbol{\omega}} \\
\dot{\boldsymbol{\omega}_s} \\
\dot{\boldsymbol{\omega}_g} \\
\dot{\q}
\end{array} \right]
& = &
\left[ \begin{array}{c}
-\J_b^{-1} \left[ \A_t  \J_{s} \boldsymbol{\Omega}_{s} \boldsymbol{\omega}_g
 + { \boldsymbol{\omega}} \times 
(\J_b \boldsymbol{\omega} +\A_s \J_s \boldsymbol{\omega}_s
+ \A_g \J_g \boldsymbol{\omega}_g)
\right]  \\
\0   \\
\0   \\
\g(\q, \boldsymbol{\omega})
\end{array} \right] +
\left[ \begin{array}{c}
\J_b^{-1} \left( \A_s \t_s  + {\A}_g \t_{g}+ \t_e \right) \\
-\J_s^{-1} \t_s  \\
-\J_g^{-1} \t_g  \\
\0
\end{array} \right] \nonumber \\
& = & 
\F ( \boldsymbol{\omega},
\boldsymbol{\omega}_g,\boldsymbol{\omega}_s, \q, t) 
+ \G (\t_s,\t_g,\t_e, t),
\label{CMGtimeVarying}
\end{eqnarray}
or simply
\begin{eqnarray}
\dot{\x} =\F(\x,\boldsymbol{\gamma}(t))
+\G(\u,\t_e, \boldsymbol{\gamma}(t)),
\end{eqnarray}
where the state variable vector is 
$\x= [ {\boldsymbol{\omega}}^{\Tr},\boldsymbol{\omega}_s^{\Tr},
\boldsymbol{\omega}_g^{\Tr}, \q^{\Tr}]^{\Tr}$, the control 
variable vector is $\u = [\t_s^{\Tr},\t_g^{\Tr}]^{\Tr}$, 
disturbance torque vector is $\t_e$,
and $\F$ and $\G$ are functions of time $t$ because the parameters 
of $\boldsymbol{\omega}$, $\boldsymbol{\omega}_s$, 
$\boldsymbol{\omega}_g$, $\q$, $\A_s$ and $\A_t$ are functions 
of time $t$. The system dimension
is $n=2N+6$. The control input dimension is $2N$.

\section{Spacecraft attitude control using variable-speed CMG}

We consider two design methods for spacecraft attitude control 
using variable-speed CMGs. But first, we approximate the 
nonlinear time-varying spacecraft control system model by a
linear time-varying spacecraft control system model near the
equilibrium point $ \boldsymbol{\omega} =\0$, 
$ \boldsymbol{\omega}_s = \0 $, 
$ \boldsymbol{\omega}_g =\0$, and $ \q = \0$ so that
an effective design considering system performance can be 
carried out using the simplified linear time-varying model. 
Denote the equilibrium by $\x_e =\0 = [\boldsymbol{\omega}^{\Tr},
\boldsymbol{\omega}_s^{\Tr},\boldsymbol{\omega}_g^{\Tr},
\q^{\Tr} ]^{\Tr} $ and
\begin{eqnarray} 
\F_1=-\J_b^{-1} \left[ \A_t  \J_{s} \boldsymbol{\Omega}_{s} \boldsymbol{\omega}_g
 + {\boldsymbol{\omega}} \times 
(\J_b \boldsymbol{\omega} +\A_s \J_s \boldsymbol{\omega}_s
+ \A_g \J_g \boldsymbol{\omega}_g)
\right], \hspace{0.1in}
\F_2=\F_3 =\0,
\hspace{0.1in} \F_4=\g(\q, \boldsymbol{\omega}),
\end{eqnarray}
\begin{equation}
\G_1=\J_b^{-1} \left(  \A_s \t_s  + {\A}_g \t_{g}+ \t_e \right), 
\hspace{0.1in}
\G_2 = - \J_s^{-1} \t_s,  \hspace{0.1in}
\G_3 = - \J_g^{-1} \t_g,  \hspace{0.1in}
\G_4 =  \0.
\end{equation}
Taking partial derivative for $\F_1$, we have 
\begin{equation}
\frac{\partial \F_1}{\partial \boldsymbol{\omega} }
= \J_b^{-1} [ (\A_s \J_s \boldsymbol{\omega}_s )^{\times}
+ (\A_g \J_g \boldsymbol{\omega}_g )^{\times}
- \boldsymbol{\omega}^{\times} \J_b 
+ (\J_b \boldsymbol{\omega} )^{\times} ]
:=\F_{11},
\label{F11}
\end{equation}
\begin{equation}
\frac{\partial \F_1}{\partial \boldsymbol{\omega}_s } 
= 
-\J_b^{-1} [ \A_t \J_s \boldsymbol{\Omega}_g
+ \boldsymbol{\omega}^{\times} \A_s \J_s ]
:=\F_{12}, 
\label{F12}
\end{equation}
\begin{equation}
\frac{\partial \F_1}{\partial \boldsymbol{\omega}_g }
= -\J_b^{-1} [ \A_t \J_s \boldsymbol{\Omega}_{s}
+ \boldsymbol{\omega}^{\times} \A_g \J_g ]
:=\F_{13},
\label{F13}
\end{equation}
\begin{equation}
\frac{\partial \F_1}{\partial \q } 
=\0.
\label{F14}
\end{equation}
Taking partial derivative for $\F_4$, we have
\begin{equation}
\frac{\partial \F_4}{\partial \boldsymbol{\omega}} 
=\frac{1}{2}  \left[  \begin{array} {ccc} 
\sqrt{1-q_1^2-q_2^2-q_3^2} & -q_3  & q_2 \\
q_3 & \sqrt{1-q_1^2-q_2^2-q_3^2} & -q_1 \\
-q_2  &  q_1  & \sqrt{1-q_1^2-q_2^2-q_3^2} \\
\end{array} \right]_{\substack{ \q \approx 0}}
\approx  \frac{1}{2} ( \I + \q^{\times})
:=\F_{41},
\label{F41}
\end{equation}
since $q_0=\sqrt{1-q_1^2-q_2^2-q_3^2}$ and
$\frac{\partial q_0}{\partial q_i}=-\frac{q_i}{q_0}$ for
$i=1,2,3$, we have
\begin{equation}
\frac{\partial \F_4}{\partial \q }
=\frac{1}{2} \left[ \begin{array}{ccc}
-\frac{q_1}{q_0} \omega_1  
& \omega_3  -\frac{q_2}{q_0} \omega_1
& -\omega_2 -\frac{q_3}{q_0} \omega_1  
\\  -\omega_3  -\frac{q_1}{q_0} \omega_2
&  -\frac{q_2}{q_0} \omega_2
&  \omega_1 -\frac{q_3}{q_0} \omega_2
\\  \omega_2 -\frac{q_1}{q_0} \omega_3
&   -\omega_1 -\frac{q_2}{q_0} \omega_3
&   -\frac{q_3}{q_0} \omega_3
\end{array} \right]_{\substack{ \boldsymbol{\omega} \approx 0 \\ \q \approx 0}} 
\approx - \frac{1}{2} \boldsymbol{\omega}^{\times}
:=\F_{44}.
\label{F44}
\end{equation}
Therefore, the linearized time-varying model is given by
\begin{eqnarray}
\left[ \begin{array}{c}
\dot{\boldsymbol{\omega}} \\
\dot{\boldsymbol{\omega}_s} \\
\dot{\boldsymbol{\omega}_g} \\
\dot{\q}
\end{array} \right]
& = &
\left[ \begin{array}{cccccc}
\F_{11}  & \F_{12}  &  \F_{13}  & \0   \\
\0 & \0 & \0 & \0   \\
\0 & \0 & \0 & \0   \\
\F_{41} & \0 & \0 & \F_{44}   \\
\end{array} \right]
\left[ \begin{array}{c}
{\boldsymbol{\omega}} \\
{\boldsymbol{\omega}_s} \\
{\boldsymbol{\omega}_g} \\
{\q}
\end{array} \right]    +
\left[ \begin{array}{cc}
\J_b^{-1} \A_s  & \J_b^{-1} {\A}_g  \\
-\J_s^{-1} & \0 \\
\0 & -\J_g^{-1} \\
\0 & \0 
\end{array} \right] 
\left[ \begin{array}{c}
{\t_s}  \\
{\t_g}  
\end{array} \right]
+ \left[ \begin{array}{c}
\J_b^{-1} \\  \0  \\ \0  \\ \0  
\end{array} \right] {\t_e}
\nonumber \\ & = & 
\A \x + \B \u +\C \t_e,
\label{CMGlinear}
\end{eqnarray}
where $\C$ is a time-invariant matrix. The linearized system 
is time-varying because ${\boldsymbol{\omega}}$, 
$\boldsymbol{\omega}_s$, $\boldsymbol{\omega}_g$, 
$\q$, $\A_s$ and $\A_t$ in $\A$ and $\B$ are all
functions of $t$. 
\begin{remark}
It is worthwhile to note that the linearized system matrices 
$\A$, $\B$, and $\C$ will be time-invariant if we approximate the
linear system at the equilibrium point of the origin ($\x_e=\0$). 
However, such a linear time invariant system will not be 
controllable. Therefore, we take the first order approximation 
for $\A$ and $\B$, which leads to a controllable linear 
time-varying system.
\end{remark}
In theory, given $\A_{s_0}$,
$\A_{t_0}$, and ${\boldsymbol{\omega}_g}$, $\A_s$ and $\A_t$ 
can be calculated by the integration of (\ref{CMGmatrixDerivative}). 
But using (\ref{diagcs}) and (\ref{ford1}) is a better method 
because it ensures that the columns of 
$\A_s$ and $\A_t$ are unit vectors as required. Notice that
the $i$th column of $\A_s$ and the $i$th column of $\A_t$, 
$i=1,\ldots,n$, must be
perpendicular to each other, an even better method to update
$\A_t$ is to use the cross product 
\begin{equation}
\t_i = \g_i \times \s_i, \hspace{0.1in} i=1,\ldots,n,
\label{crossperpendicular}
\end{equation}
to prevent $\t_i$ and $\s_i$ from being losing perpendicularity 
due to the numerical error accumulation. In simulation, integration of (\ref{dotGamma}) can be used to obtain 
$\boldsymbol{\gamma}$ which is needed in the computation of 
(\ref{diagcs}), but in engineering practice, the encoder 
measurement should be used to get $\boldsymbol{\gamma}$.

Assuming that the closed-loop linear time-varying system   
is given by
\begin{equation}
\dot{\x} = \bar{\A}(t) \x(t), \hspace{0.1in} \x(t_0) =\x_0.
\label{timeVarying}
\end{equation}
It is well-known that even if all the eigenvalues of 
$\bar{\A}(t)$, denoted by $\Re [\lambda (t)]$, are in the left 
half complex plane for all $t$, the system may not be stable 
\cite[pages 113-114]{rugh93}. But the following
theorem (cf. \cite[pages 117-119]{rugh93}) provides a nice 
stability criterion for the closed-loop system 
(\ref{timeVarying}). 
\begin{theorem}
Suppose for the linear time-varying system (\ref{timeVarying})
with $\bar{\A}(t)$ continuously differentiable there exist finite 
positive constants $\alpha$, $\mu$ such that, for all 
$t$, $\| \bar{\A}(t) \| \le \alpha $ and every point-wise eigenvalue 
of $\bar{\A}(t)$ satisfies $\Re [\lambda (t)] \le -\mu$.
Then there exists a positive constant $\beta$ such that if the
time derivative of $\bar{\A}(t)$ satisfies 
$\| \dot{\bar{\A}}(t) \| \le \beta$ for all $t$, the 
state equation is uniformly exponentially stable.
\label{rughTimeVarying}
\end{theorem}

This theorem is the theoretical base for the linear time-varying
control system design. We need at least that 
$\Re [\lambda (t)] \le -\mu$ holds.

\subsection{Gail scheduling control}

Gain scheduling control design is fully discussed in 
\cite{rugh90} and it seems to be applicable to this LTV system. 
The main idea of gain scheduling is: 1) select a set of fixed
parameters' values, which represent the range of the plant
dynamics, and design a linear time-invariant gain for each;
and 2) in between operating points, the gain is interpolated 
using the designs for the fixed parameters' values that cover
the operating points. As an example,
for $i=1,\ldots,N$, let $\gamma_{i} \in \{ 2\pi/p_{\gamma}, 
4\pi/p_{\gamma}, \cdots,2\pi \}$ be a set of $p_{\gamma}$ fixed 
points equally spread in $[0,2\pi]$. Then, for $N$ CMGs, there 
are $p_{\gamma}^N$ possible fixed parameters' combinations. 
For example, if $N=4$ and $p_{\gamma}=8$, we can represent 
the grid composed of these fixed points in a matrix form 
as follows:
\begin{equation}
\left[ \begin{array}{cccccccc}
\pi/4 & \pi/2 & 3\pi/4 & \pi & 5\pi/4 & 3\pi/2 & 7\pi/4 & 2\pi \\
\pi/4 & \pi/2 & 3\pi/4 & \pi & 5\pi/4 & 3\pi/2 & 7\pi/4 & 2\pi \\
\pi/4 & \pi/2 & 3\pi/4 & \pi & 5\pi/4 & 3\pi/2 & 7\pi/4 & 2\pi \\
\pi/4 & \pi/2 & 3\pi/4 & \pi & 5\pi/4 & 3\pi/2 & 7\pi/4 & 2\pi 
\end{array} \right],
\label{grid}
\end{equation}
and each fixed $\boldsymbol{\gamma}$ is a vector composed of 
$\gamma_i$ ($i=1,2,3,4$) which can be any element of $i$th row.
If $\boldsymbol{\gamma}$ is not a fixed point, we have 
$\gamma_i \in [\kappa(i), \kappa(i) +1]$ for all 
$i \in [1,\cdots,N-1]$. Assume that $\gamma_i$ is in 
the interior of $(\kappa(i), \kappa(i) +1)$ for all
$i \in [1,\cdots,N-1]$. Then, $\boldsymbol{\gamma}$ meets the
following conditions:
\begin{equation}
\boldsymbol{\gamma} =
\left[ \begin{array}{c}
\gamma_1 \in (\kappa(1), \kappa(1)+1) \\ 
\vdots \\ 
\gamma_N  \in (\kappa(N), \kappa(N)+1)
\end{array} \right].
\label{vertex}
\end{equation}
Using the previous example of (\ref{grid}), if $\boldsymbol{\gamma}
=\left[ \frac{5\pi}{8},\frac{3\pi}{8},\frac{7\pi}{16},
\frac{15\pi}{8} \right]^{\Tr}$, then
$\boldsymbol{\gamma} \in \left[ (\frac{\pi}{2},\frac{3\pi}{4}),
(\frac{\pi}{4},\frac{\pi}{2}),
(\frac{\pi}{4},\frac{\pi}{2}),
(\frac{7\pi}{4},{2\pi})  \right]^{\Tr}$.
To use gain scheduling control, we need also to consider fixed 
points for ${\boldsymbol{\omega}}$, 
$\boldsymbol{\omega}_s$, $\boldsymbol{\omega}_g$, and
$\q$ in their possible operational ranges. 
Let $p_w$, $p_{w_s}$, $p_{w_g}$, and $p_q$ be the
number of the fixed points for ${\boldsymbol{\omega}}$, 
$\boldsymbol{\omega}_s$, $\boldsymbol{\omega}_g$, and
$\q$. The total vertices for the entire polytope (including a
grid of all possible time-varying parameters) will be
$p_{\gamma}^N p_w^3 p_{w_s}^N p_{w_g}^N p_q^3$. 

For each of these ($p_{\gamma}^N p_w^3 p_{w_s}^N p_{w_g}^N p_q^3$) fixed models, we need conduct a control design to 
calculate the feedback gain matrix for the ``frozen'' model. 
If the system (\ref{CMGlinear}) at time $t$ happens to have 
all parameters equal to 
the fixed points, we can use a ``frozen'' feedback gain to control 
the system (\ref{CMGlinear}). Otherwise, we need to construct
a gain matrix based on $2^{3N+6}$  ``frozen'' gain matrices.
Assuming that each parameter has some moderate number of fixed 
points, say $8$, and the control system has $N=4$ gimbals, the 
total number of the fixed models will be $8^{18}$, each needs
to compute a feedback matrix, an impossibly computational task.

\subsection{Model Predictive Control}

Unlike the gain scheduling control design in which most computation
is done off-line, model predictive control computes the feedback 
gain matrix on-line for the linear system (\ref{CMGlinear})
in which $\A$ and $\B$ matrices are updated in every 
sampling period. It is straightforward to verify that for any 
given $\boldsymbol{\gamma}$, if $\x \neq \x_e$, the linear
system (\ref{CMGlinear}) is controllable. In theory, one can use 
either robust pole assignment \cite{yt93,tity96}, or LQR design 
\cite{lvs12}, or $\H_{\infty}$ design \cite{zdg96} for the
on-line design, but $\H_{\infty}$ design costs significant more
computational time and should not be considered for this on-line 
design problem. Since LTV system design should meet the condition 
of $\Re [\lambda (t)] \le -\mu$ required in Theorem 
\ref{rughTimeVarying}, robust pole assignment design is 
clearly a better choice than LQR design for this purpose. 
Another attrictive feature of the robust pole assignment design
is that the perturbation of the closed loop eigenvalues between 
sampling period are expected to be small.
It is worthwhile to note that a robust pole assignment design 
\cite{tity96} minimizes an upper bound of $\H_{\infty}$ norm which
means that the design is robust to the modeling error and reduces 
the impact of disturbance torques on the system output 
\cite{yang96,yang14}. Additional merits about this method, such
as computational speed which is important for the on-line design, 
is discussed in \cite{psnyst14}. Therefore, we use the method of \cite{tity96} in the proposed design.

The proposed design algorithm is given as follows:

\begin{algorithm} {\ } \\ 
Data: $\J_b$, $\J_s$, $\J_g$, and $\A_g$.  \hspace{0.1in} {\ } \\
Initial condition: $\x=\x_0$, $\boldsymbol{\gamma}=\boldsymbol{\gamma}_0$, $\A_{s_0}$, 
and $\A_{t_0}$.   {\ } \\
\begin{itemize}
\item[] Step 1: Update $\A$ and $\B$ based on the latest
$\boldsymbol{\gamma}$ and $\x$.
\item[] Step 2: Calculate the gain $\K$ using robust pole assignment
algorithm {\tt robpole} (cf. \cite{tity96}).
\item[] Step 3: Apply feedback $\u=\K \x$ to (\ref{CMGtimeVarying}) 
or (\ref{CMGlinear}).
\item[] Step 4: Update $\boldsymbol{\gamma}$ and 
$\x = [ {\boldsymbol{\omega}}^{\Tr},
\boldsymbol{\omega}_s^{\Tr},
\boldsymbol{\omega}_g^{\Tr}, \q^{\Tr}]^{\Tr}$.
Go back to Step 1.
\end{itemize}
\label{onLine}
\end{algorithm}

\section{Simulation test}

%

\begin{figure}[ht]
\centering
\includegraphics[height=8cm,width=8cm]{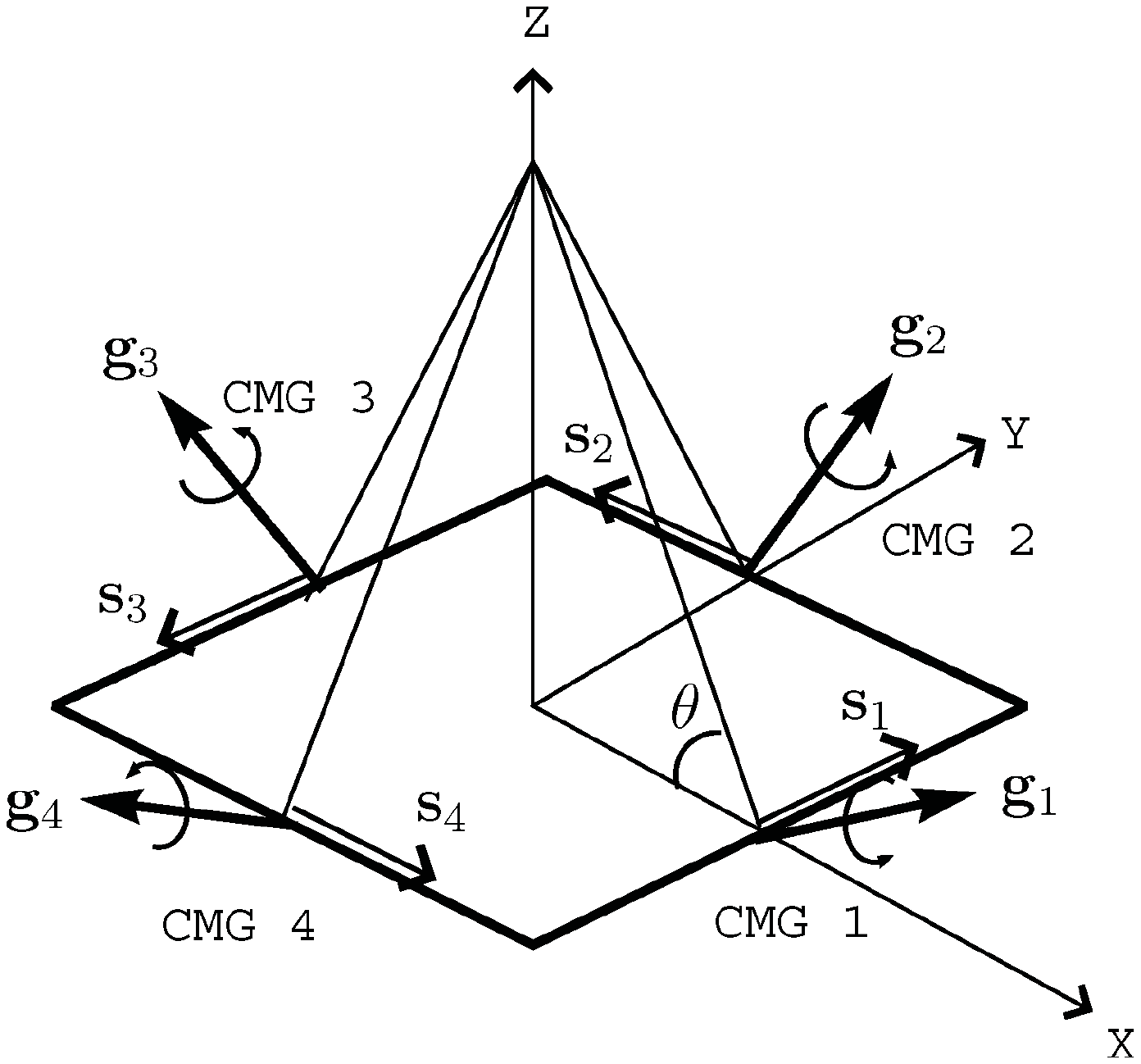}
\caption{VSCMG system with pyramid configuration concept.}
\label{fig:pyramid1}
\end{figure}

\begin{figure}[ht]
\centering
\includegraphics[height=8cm,width=8cm]{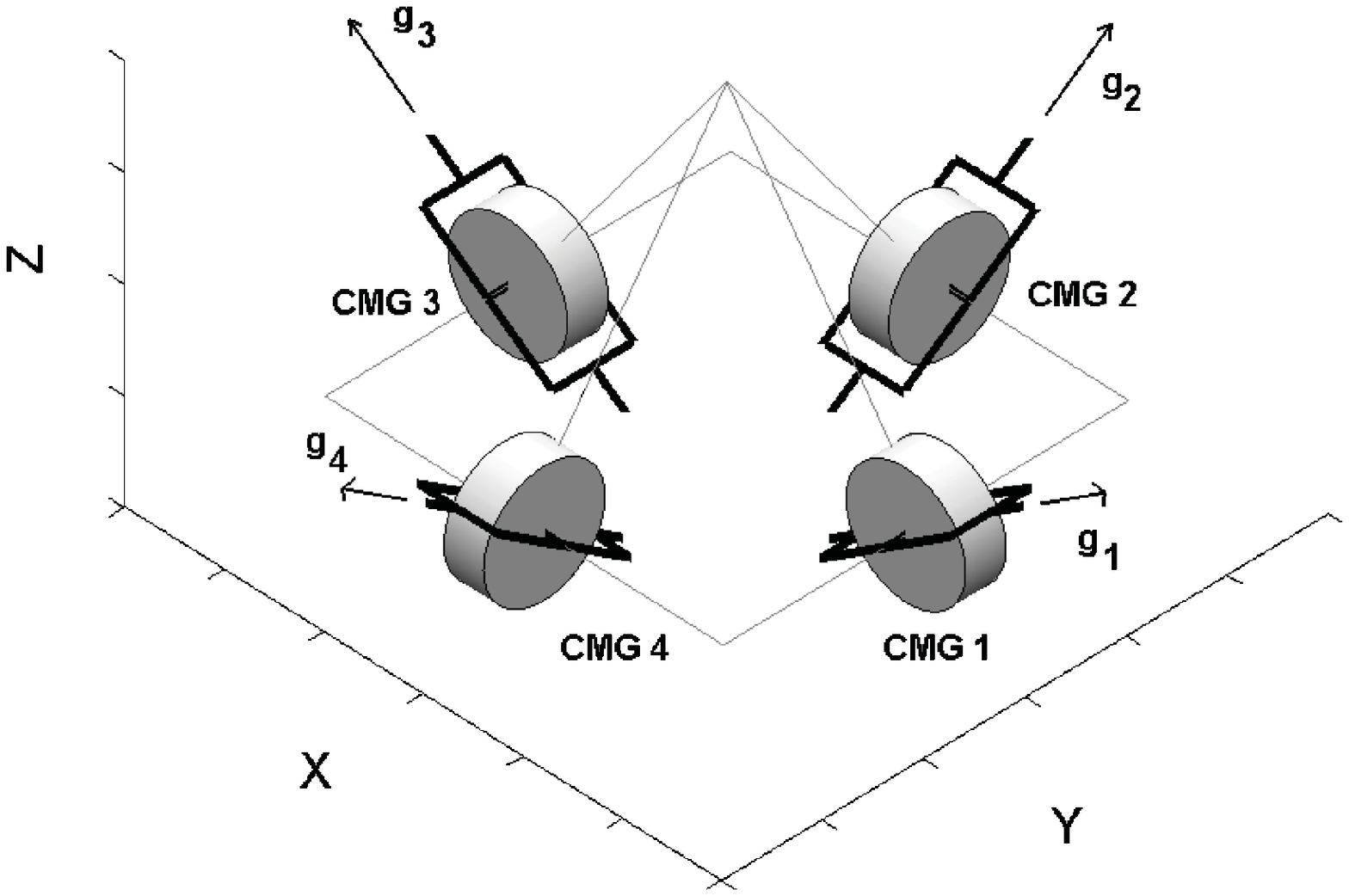}
\caption{VSCMG system with pyramid configuration.}
\label{fig:pyramid2}
\end{figure}

\begin{figure}[ht]
\centering
\includegraphics[height=6cm,width=8cm]{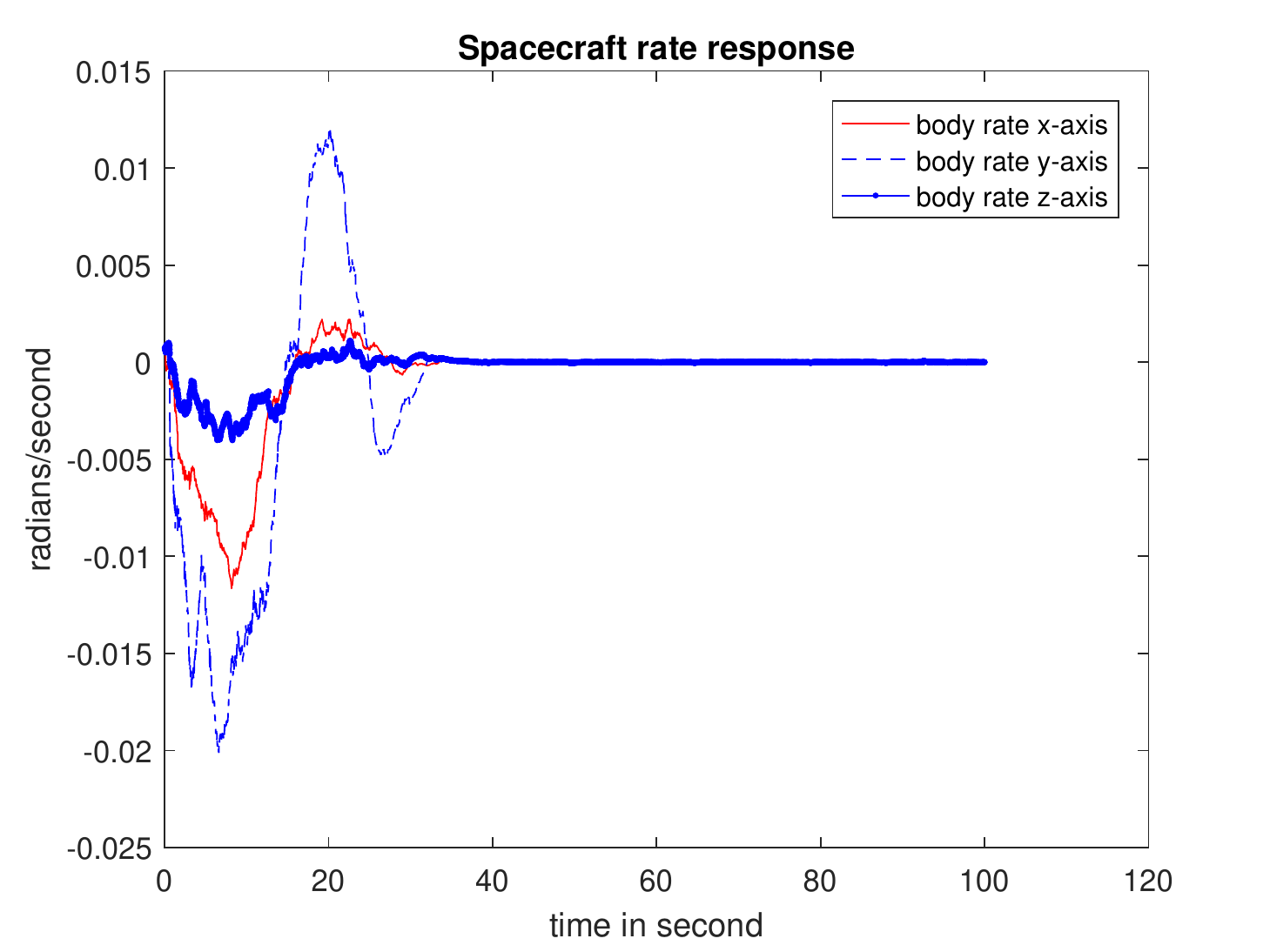}
\caption{Spacecraft body rate response.}
\label{fig:bodyRate}
\end{figure}

\begin{figure}[ht]
\centering
\includegraphics[height=6cm,width=8cm]{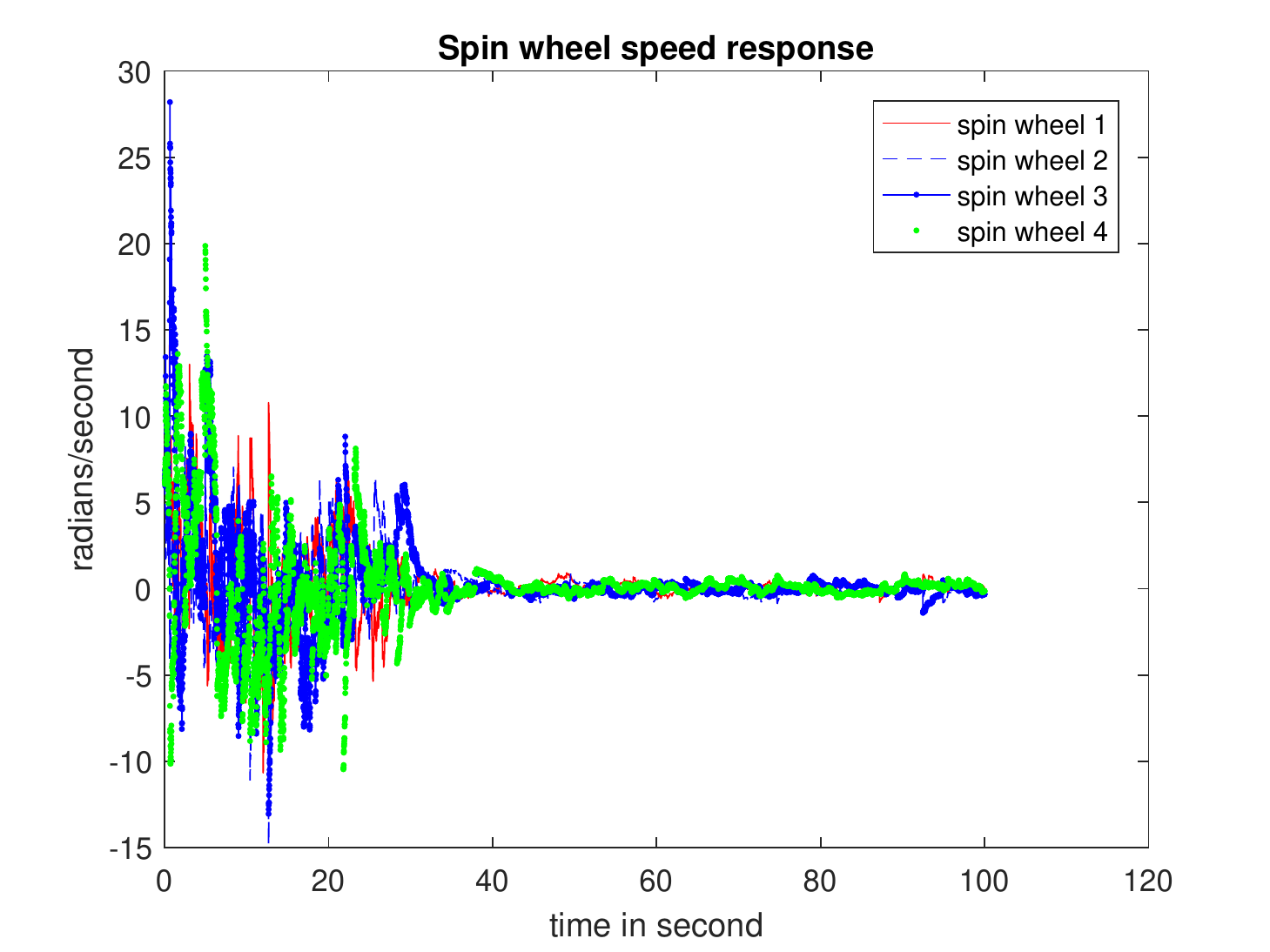}
\caption{Spin wheel response.}
\label{fig:spinWheel}
\end{figure}

\begin{figure}[ht]
\centering
\includegraphics[height=6cm,width=8cm]{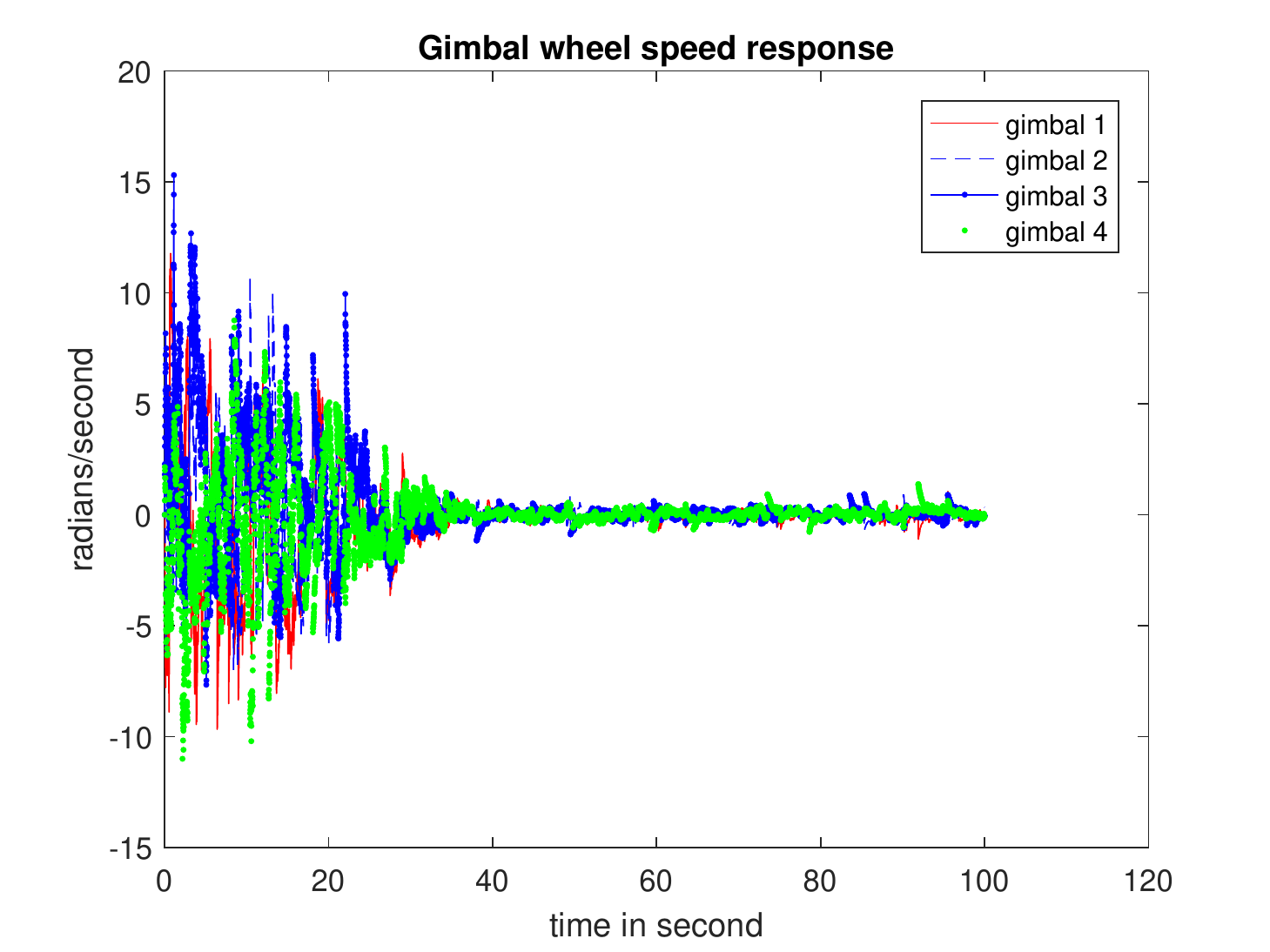}
\caption{Gimbal wheel response.}
\label{fig:gimbalWheel}
\end{figure}

\begin{figure}[ht]
\centering
\includegraphics[height=6cm,width=8cm]{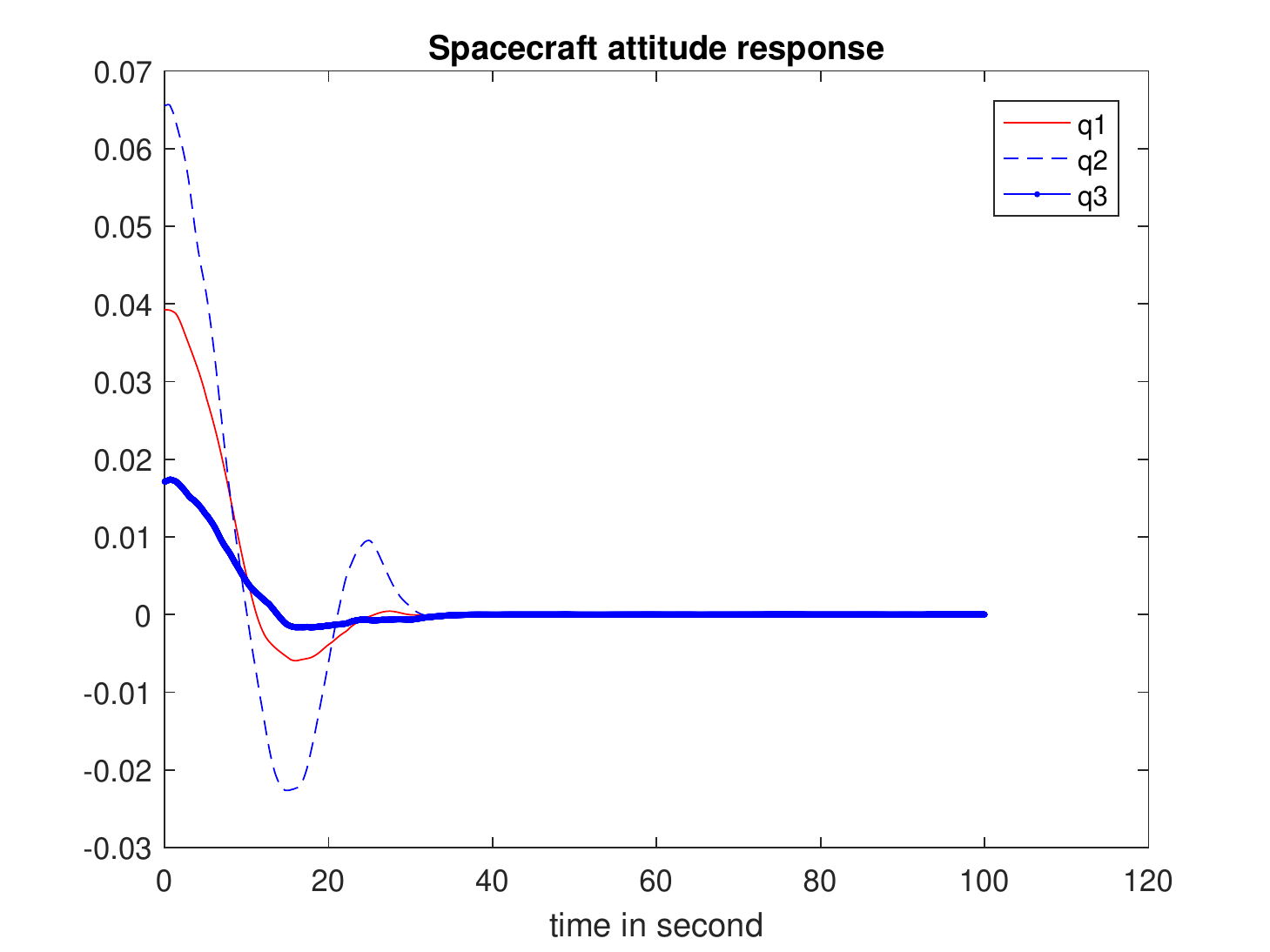}
\caption{Reaction wheel response $\Omega_1$, $\Omega_2$, and $\Omega_3$.}
\label{fig:quaternionRes}
\end{figure}

The proposed design method is simulated using the data in 
\cite{jt04,yt04,yt02}. We assume that the four variable-speed 
CMGs are mounted in pyramid configuration as shown in Figures
\ref{fig:pyramid1} and \ref{fig:pyramid2}.
The angle of each pyramid side to its base is $\theta=54.75$ 
degree; the inertia matrix of the spacecraft is given by 
\cite{yt02} as
\begin{equation}
\J_b = \left[ \begin{array}{ccc}
15053 & 3000 & -1000 \\
3000 & 6510 & 2000 \\
-1000 & 2000 & 11122
\end{array} \right] 
\hspace{0.1in} \mbox{kg}\cdot m^2.
\end{equation}
The spin axis inertial matrix is given by
$\J_s = \diag ( 0.7, 0.7, 0.7, 0.7) \hspace{0.04in} \mbox{kg}\cdot m^2$ 
and the gimbal axis inertia matrix is given by 
$\J_g=\diag( 0.1, 0.1, 0.1, 0.1) \hspace{0.04in} \mbox{kg}\cdot m^2$.
The initial wheel speeds are $2 \pi$ radians per second for all
wheels. The initial gimbal speeds are all zeros. The initial spacecraft body rate vector is randomly generated by Matlab 
$rand(3,1) *10^{-3}$ and the initial spacecraft attitude vector 
is a reduced quaternion randomly generated by Matlab 
$rand(3,1) *10^{-1}$. The gimbal axis matrix is fixed and given
by \cite{yt04} (cf. Figures \ref{fig:pyramid1} and \ref{fig:pyramid2}.)
\begin{equation}
\A_g = \left[ \begin{array}{cccc}
\sin(\theta) & 0 & -\sin(\theta) & 0 \\
0 & \sin(\theta) & 0 & -\sin(\theta) \\
\cos(\theta) & \cos(\theta) & \cos(\theta) & \cos(\theta) 
\end{array} \right]
\label{gimbalM}
\end{equation}
The initial wheel axis matrix can be obtained using Figures 
\ref{fig:pyramid1} and \ref{fig:pyramid2} and is given by
\begin{equation}
\A_s = \left[ \begin{array}{cccc}
0 & -1 & 0 & 1 \\
1 & 0 & -1 & 0 \\
0 & 0 & 0 & 0
\end{array} \right]
\label{gimbalW}
\end{equation}
The initial transverse matrix $\A_t$ can be obtained by the 
method of (\ref{crossperpendicular}). The desired or designed
closed-loop poles are selected as $\{ -0.2 -0.8, -0.2 \pm 0.1i, -0.6\pm 0.1 i,-1.5 \pm i, -1.6 \pm i,-1.7 \pm i,-1.8 \pm i\}$. The simulation test results are given 
in Figures \ref{fig:bodyRate}-\ref{fig:quaternionRes}. Clearly,
the designed controller stabilizes the system with good 
performance.

%
%
%

\section{Conclusions} 
 
In this paper, we proposed a new operational concept for 
variable-speed CMGs. This new concept allows us to simplify
the nonlinear model of the spacecraft  attitude control using 
variable-speed CMGs to a linear time-varying model. Although
this LTV model is significantly simpler than the original nonlinear
model, there are still many time-varying parameters in the 
simplified model. Two LTV control system design methods,
the gain scheduling design and model predictive control design,
are investigated. The analysis shows that model predictive 
control is better suited for spacecraft control using variable-speed
CMGs. An efficient robust pole assignment algorithm is used in
the on-line feedback gain matrix computation. Simulation test
demonstrated the effectiveness of the new concept and control 
system design method.


\end{document}